\documentclass{article}
\usepackage[all]{xy}
\usepackage[latin1]{inputenc}
\usepackage{amscd}
\usepackage{graphicx}
\usepackage{rotating}
\usepackage[english]{babel}
\usepackage{amsmath,amssymb,amsfonts}
\usepackage{setspace}                    
\usepackage{latexsym}
\usepackage{color}
\usepackage{epsfig}
\begin{document}

\title{The genus of the configuration spaces for Artin groups of affine type}
\author{D. Moroni \thanks{Institute of Information Science and Technologies
(ISTI), National Research Council of Italy (CNR), Pisa}, M. Salvetti \thanks
{Corresponding author: Mario Salvetti, Department of  Matemathics, 
University of Pisa, e-mail: salvetti@dm.unipi.it}, A. Villa \thanks {Department of Mathematical Sciences , ISI Foundation, e-mail: andrea.villa@isi.it 
}}
\date {}
\maketitle

\begin{abstract}
Let $(\mathbf W,S)$ be a Coxeter system, $S$ finite, and let $\mathbf G_{\mathbf W}$ be the associated Artin group.
One has {\it configuration spaces} $\mathbf Y,\ \mathbf Y_{\mathbf W},$ where $\mathbf G_{\mathbf W}=\pi_1(\mathbf Y_{\mathbf W}),$ and a natural   $\mathbf W$-covering $f_{\mathbf W}:\ \mathbf Y\to\mathbf Y_{\mathbf W}.$ The {\it Schwarz genus}  $g(f_{\mathbf W})$ is a natural topological invariant to consider. In \cite{salvdec2} it was computed for all finite-type Artin groups, with the exception of case $A_n$ (for which see \cite{vassiliev},\cite{salvdecproc3}). 
In this paper we generalize this result by computing the Schwarz genus for a class of Artin groups, which includes the affine-type Artin groups.
Let $K=K(\mathbf W,S)$ be the simplicial scheme of all subsets  $J\subset S$ such that the parabolic group $\mathbf W_J$ is finite.  We introduce the class of groups  for which $dim(K)$ equals the homological dimension of $K,$ and  we show that   $g(f_{\mathbf W})$ is always the maximum possible for such class of groups. For affine Artin groups, such maximum reduces to the rank of the group. In general, it is given by 
$dim(\mathbf X_{\mathbf W})+1,$  where $\mathbf X_{\mathbf W}\subset \mathbf Y_{\mathbf W}$ is   a well-known $CW$-complex which has the same homotopy type as $\mathbf Y_{\mathbf W}.$
\end{abstract}

\newcommand{\oC}[0]{\overline{C}}
\newcommand{\oF}[0]{\overline{F}}
\newcommand{\oG}[0]{\overline{G}}
\newcommand{\Y}[0]{\mathbf Y^{(d)}}
\newcommand{\arr}[0]{{\mathcal A}}
\newcommand\X{\mathbf X}
\newcommand\Xd{\mathbf X^{(d)}}
\newcommand\Xdd{\mathbf X^{(d+1)}}
\newcommand\Sd{\mathbf S^{(d)}}
\newcommand\Sdd{\mathbf S^{(d+1)}}
\renewcommand{\S}[0]{\mathbf{S}}
\newcommand\W{\mathbf W}
\newcommand\WG{\mathbf W_{\Gamma}}
\newcommand\WS{\mathbf W^{\Gamma}}
\newcommand\XW{\mathbf X_{\mathbf W}}
\newcommand\Fi{\mathbf \varPhi}
\renewcommand\>{\succ}
\newcommand\po{\lessdot}
\newcommand{\poi}[1]{\lessdot_{#1}^{op}}
\newcommand\Q{\mathbf Q}
\newcommand\Qdd{\mathbf Q^{(d+1)}}
\newcommand\Qd{\mathbf Q^{(d)}}
\newcommand\Qk{\mathbf Q^{(k)}}
\newcommand\join{\bigcup}
\newcommand{\strat}[1]{\Fi}
\newcommand\strk{\Fi ^{(k)}}
\newcommand\strd{\Fi ^{(d)}}
\newcommand\str{\Fi^{(d+1)}}
\newcommand\strY{\Fi^{(d+1)}_0}
\newcommand\strA{\Fi^{(d+1)}_{+}}
\newcommand\strR{\Fi^{(d)}}
\newcommand\prR{pr_{\Re}}
\newcommand\s{\hat}
\newcommand{\Fp}[0]{\Fc^{\prime}}
\newcommand\Fs{\Fc^{\prime\prime}}
\newcommand\cam{\mathbf \varPhi_0}
\newcommand\g{\gamma}
\newcommand\gp{\gamma ^{\prime}}
\newcommand\G{\Gamma}
\newcommand\Gp{\Gamma ^{\prime}}
\newcommand\FO{\{\mathbf 0\}}
\newcommand\GF{ \varGamma}
\newcommand\Gf{\varGamma^{\prime}}
\newcommand\pW{\pi_{\W}}
\newcommand\gk{\hat\g}
\newcommand\GW{\mbf{G_W}}
\newcommand\PW{\mbf{PG_W}}

\newcommand{\A}[0]{\mathcal{A}}
\newcommand{\M}[0]{\mathcal{M}}
\newcommand{\Md}[1]{\mathcal{M}^{(#1)}}
\newcommand{\CS}[2]{\mathcal{#1}^{(#2)}}
\newcommand{\La}[0]{\mathcal{L}(\A)}
\renewcommand{\L}[0]{\mathcal{L}}
\newcommand{\F}[0]{\mathbb{F}}
\newcommand\Fc{\mathcal{F}}
\newcommand\Gc{\mathcal{G}}
\newcommand{\Z}[0]{\mathbb{Z}}
\newcommand{\R}[0]{\mathbb{R}}
\newcommand{\C}[0]{\mathbb{C}}
\newcommand{\CW}[0]{\mathcal{C}}
\newcommand{\N}[0]{\mathbb{N}}
\newcommand{\K}[0]{\mathbb{K}}
\newcommand{\St}[0]{\mathcal S}
\newcommand{\Tbf}[0]{\mathbf{T}}
\newcommand{\Cal}[1]{\mathcal{#1}}
\newcommand{\vs}[0]{\vspace{3mm}}
\newcommand{\ph}[0]{\varphi}
\newcommand{\la}[0]{\lambda}
\newcommand{\bd}[1]{\textbf{#1}}
\newcommand{\e}{\bd e}
\renewcommand{\lessdot}{\vartriangleleft}
\newcommand{\MS}{\tilde{\bd S}}
\newcommand{\te}{\theta}
\newcommand{\Te}{\mathbf{\theta}}
\newcommand{\q}[1]{\mbox{\bfseries{\textit{#1}}}}
\newcommand{\itl}[1]{\textit{#1}}
\newcommand{\p}[0]{p}
\newcommand{\geets}[0]{\longleftarrow}
\newcommand{\too}[0]{\longrightarrow}
\newcommand{\sst}[0]{\subset}
\newcommand{\cl}[1]{\mathcal{#1}}
\newcommand{\into}[0]{\hookrightarrow}
\newcommand{\codim}[0]{\mbox{codim }}
\newcommand{\diag}[0]{\Sigma}
\newcommand{\ea}[0]{\underline{\epsilon}}
\newcommand{\eb}[0]{\overline{\epsilon}}
\newcommand{\im}[0]{\mbox{im }}
\newcommand{\az}[0]{\"{a}}
\newcommand{\ls}[0]{\mathcal L}
\newcommand{\ul}[0]{\underline}
\newcommand{\bin}[2]{  \left(  \begin{array}{c}  #1 \\ #2  \end{array}
  \right)  }
\newcommand{\qbin}[2]{  \left[  \begin{array}{c}  #1 \\ #2
    \end{array}  \right]  }
\newcommand{\qed}[0]{\hspace{\stretch{1}}$\Box$}
\newcommand{\eq}[1][r]
       {\ar@<-3pt>@{-}[#1]
        \ar@<-1pt>@{}[#1]|<{}="gauche"
        \ar@<+0pt>@{}[#1]|-{}="milieu"
        \ar@<+1pt>@{}[#1]|>{}="droite"
        \ar@/^2pt/@{-}"gauche";"milieu"
        \ar@/_2pt/@{-}"milieu";"droite"}
\newcommand{\imm}[1][r] {\ar@{^{(}->}[#1]}
\newcommand{\D}[0]{D^{\scriptscriptstyle{(0)}}}
\newcommand{\dd}[0]{d^{\scriptscriptstyle{(0)}}}
\renewcommand{\ni}[0]{\noindent}
\newcommand{\Ga}[0]{\Gamma}

\newcommand{\flo}[1]{\lfloor #1 \rfloor}
\newcommand{\pr}{\mathrm{Prod}}
\providecommand{\poi}[1]{{\mathbb{W}}_{#1}}
\providecommand{\qed}{{\flushright $\square$\\}}
\renewcommand{\ll}{{ \mathcal{L}_2 (\mathbb{R}^3, d \mu)}}
\providecommand{\p}{{\mathbf{p}}}
\renewcommand{\L}{{\mathbf{L}}}
\renewcommand{\i}{{\mathrm{i}}}
\providecommand{\Br}{{\mathrm{Br}}}
\providecommand{\A}{{\mathbf{A}}}
\newcommand{\lie}[2]{{[ #1 , #2  ]}}
\newcommand{\bra}[2]{{\langle #1, #2 \rangle}}
\newcommand{\no}[1]{{| #1 |}}
\newcommand{\til}[1]{\widetilde{#1}}
\newcommand{\qa}[1]{[#1]}
\newcommand{\fhi}[1]{\{ #1 \}}
\renewcommand{\ln}[1]{\bar{#1}}
\newcommand{\bo}[1]{{\mathcal{#1}}}
\newcommand{\mbf}[1]{{\mathbf{#1}}}

\newcommand{\re}{{\mathbf{R}}}
\newcommand{\B}{{\mathcal{B}}}
\newcommand{\card}[1]{|#1|}
\newcommand{\oth}{\mathrm{otherwise}}
\newcommand{\semidir}{\rtimes}
\newcommand{\de}{\partial}
\newcommand{\di}{{\mathrm{d}}}
\renewcommand{\mod}[1]{\,\mathrm{mod}(#1)}
\newcommand{\compl}{\Xi}
\newcommand{\wI}[0]{\widetilde{I}}
\newcommand{\wqbin}[2]{ \widetilde{\left[ \begin{array}{c} #1 \\ #2
    \end{array} \right] }}
\newcommand{\tss}[0]{\supseteq}
\newcommand{\pmu}[0]{{\pm 1}}
\newcommand{\CA}[0]{\overline{C}}


\newcommand{\br}{\mathbf{G_W}}
\newcommand{\pbr}[1]{PG_{#1}}
\renewcommand{\t}[1]{\tilde{#1}}
\newcommand{\rank}{\mathrm{rk}}
\newcommand{\ind}{\mathrm{Ind}}
\newcommand{\coind}{\mathrm{CoInd}}
\newcommand{\pp}[1]{{(#1)}}


\newtheorem{df}{Definition}[section]
\newtheorem{teo}{Theorem}
\newtheorem{claim}{Claim}
\newtheorem{prop}[df]{Proposition}
\newtheorem{lem}[df]{Lemma}
\newtheorem{cor}[df]{Corollary}
\newtheorem{rmk}[df]{Remark}
\newtheorem{notat}[df]{Notation}
\newtheorem{ex}{Example}
\newtheorem{pf}{Proof}
\newtheorem{for}{}
\newtheorem{conj}{Conjecture}

\newenvironment{es}[1][Example.]{\begin{trivlist}
     \item[\hskip \labelsep {\bfseries #1}]}{\end{trivlist}}
\newenvironment{dm}[1][Proof.]{\begin{trivlist} \item[\hskip
    \labelsep {\bfseries #1}]}{\end{trivlist}}
\newenvironment{dmof}[1][Proof]{\begin{trivlist} \item[\hskip
    \labelsep {\bfseries #1}]}{\end{trivlist}}
\newenvironment{os}[1][Remark.]{\begin{trivlist}
     \item[\hskip \labelsep {\bfseries #1}]}{\end{trivlist}}
\newenvironment{grafi}

\section{Introduction} 


To any   Coxeter system $(\W,S),$  $S$ finite,  one can naturally associate (see section \ref{sec2}):
\begin{enumerate}
\item[-]  a  space
$\mbf Y$ with a natural free action of $\W,$ and an {\it orbit space} $\mbf Y_{\W}$  such that  
the projection onto the quotient \ $f_{\W}:\mbf Y  \to  \mbf Y_{\W}$ \ 
is a regular covering with group $\W;$  
\item[-] an explicitly constructed  $CW$-complex $\X\subset \mbf Y$ which is a deformation retract of $\mbf Y,$ whose cells are permuted under the action of $\W,$ and a {\it finite}  orbit $CW$-complex $\X_{\W}\subset \mbf Y_{\W}$  which is a deformation retract of the orbit space, such that $f_{\W}$ restricts to a regular $\W$-covering $f_{\W}:\X\to \XW.$   
\end{enumerate}

Let $K=K(\W)$ be the simplical scheme with vertex set $S$ and simplices all $\Gamma \subset S$ such that the parabolic subgroup $\W_{\Gamma}$ is finite. Then the complex $\X$ is union of finitely many polyhedra, one for each maximal simplex in $K,$ and $\XW$ is obtained by explicit identifications on the faces of these polyhedra (see 
\cite{salvettiArtin}).  

The fundamental group of $\mbf Y_{\W}$  (resp. of $\mbf Y$)  is the {\it Artin group} \  $\GW$ of type $\W$ (resp. the {\it pure} Artin group $\PW$).  
The quotient  $\GW/ \PW$ is isomorphic to $\W.$  For example,  if the group $\W$ is the symmetric group $\Sigma_n$ with set of Coxeter generators $S:=\{(i,i+1): \ i=1,\dots,n-1\}$ then $\GW$ is the {\it braid group} $\Br_{n}$ and $\PW$ is the pure braid group $P_n.$ \ The spaces $\mbf Y_{\W}$ and  $\mbf Y$  generalize the classical {\it configuration spaces} for the braid group and the pure braid group respectively, so they can be called the {\it configuration spaces} for  $\GW$  and  for  $\PW$  respectively.

The  {\it Schwarz genus}  $g(f)$ of a locally trivial fibration $f:Y\to X$ is
the minimum cardinality of an open covering $\mathcal U$ of $X$ such that there exists a section 
of $f$ over each open set  $U\in \mathcal U$ (\cite{schwarz}; see section \ref{genus} below). The Schwarz genus $g(f_{\W})$ of the above covering can be considered  as a natural topological invariant of the Artin group. 

The main result of this paper  is the computation of $g(f_{\W})$ for a class of Artin groups including the {\it affine type} groups (i.e., 
when $(\W,S)$ is an affine Coxeter system).

For any finite Coxeter group ${\W}$ of rank $n$, it was shown in
\cite{salvdec2} that the Schwarz genus $g(f_{\W})$ reaches the
upper bound $n+1$ (coming from an easy argument in obstruction theory)
except when ${\W}=A_n$ and $n+1$ is not a prime power (see \cite{vassiliev}).
For the first non prime power $n+1=6$, in \cite{salvdecproc3} it is shown that  $g(f_{A_5})=5,$  lower than $n+1.$
Next, in \cite{arone} it is shown by different methods that $g(f_{A_n})< n+1$
when $n+1\neq p^k, 2p^k$ for $p$ prime (so, the case $2p^k$ is still open, as 
well as the precise value of the genus).

We introduce the class of {\it affine-like} Artin groups, as those groups such that the homological dimension of $K(\W)$ equals the topological dimension of $K(\W).$ For affine type Artin groups, $K(\W)$ turns out to be a sphere, so this condition is fulfilled.  For this class of groups,  we show that the genus coincides again with 
$dim(\XW)+1.$ For an affine type Artin group of rank $n+1,$   $dim(\XW)=n$ (but for other 
affine-like groups, the rank is bigger than $dim(\XW)+1$).

Notice how the above condition on $K$ contrasts  what happens for the finite type cases, where $K$ is contractible.

\section{General pictures}\label{sec2}

\subsection{Topological constructions for Artin groups}\label{subsec:topconArt}
We will consider a finitely generated   {\it Coxeter system} $(\W,S)$ \ ($S$
finite), so
\begin{equation}\label{presentationW}
\W\ =\ <s\in S\ |\ (ss')^{m(s,s')}=1>
\end{equation}
where $m(s,s')\in\N\cup\{\infty\},\ m(s,s')=m(s',s),\ m(s,s)=1$ 
(all results in this part which are not explicitly referred to the literature are taken from
\cite{bourbaki}, \cite{humphreys}, \cite{vinberg}).
We recall some general pictures. 
 
The group $\W$ can be realized as a group generated by (in general, non-orthogonal)
reflections in $\R^n,\ n=|S|.$  
Let $\A$ be the {\it reflection arrangement}, i.e.
$$\A\ =\ \{H\subset\R^n\ |\ H\ \mbox {is fixed by some reflection in}\ \W\} .$$
Consider also the stratification into {\it facets} $\strat\ := \{ F\}$  of
$\R^n$ induced by $\A$. The codimension-$0$ facets, which are the connected 
components of the complement to the arrangement,
are called {\it chambers}.
All the chambers are simplicial cones, the group acting transitively over the
set of all them.  The Coxeter generator set $S$ corresponds to the set of
reflections with respect to the {\it walls} of a fixed base-chamber $C_0.$

Let \ $U:= \W.\overline{C}_0$ \ be the orbit of the closure of the base chamber.
$U$ is called the {\it Tits cone} of the Coxeter system. 

Notice that  the closure of the chamber $\overline{C}_0$ is endowed with a
natural stratification into facets (which are still relatively open cones with
vertex $0$). When $C_0$ is the standard positive octant, each facet is given by
imposing some coordinates equal to $0,$ and the remaining coordinates positive.
 
Each reflection in $\W$ is conjugated to a reflection with respect to a wall of
$C_0.$ So,
the arrangement $\A$  consists of the orbits of the walls of
$C_0.$ Each chamber contained inside $U$ is of
the shape $w.C_0$ for a unique $w\in \W.$   Of course,  $\A$  is not locally
finite if $\W$ is infinite  (e.g. $0$ is contained in all the hyperplanes). The orbits of the facets
of $C_0$ give a ``stratification'' of $U$ into relatively open cells, also
called facets (in general, $U$ is neither open nor closed in $\R^n$).

Recall also:
\begin{enumerate}
\item $U$ is a convex cone in $\R^n$ with vertex $0.$
\item $U=\R^n$ \ iff \ $\W$ is finite
\item The stabilizer of a facet $F$ in $U$ is the subgroup $\mbf{W}_F$ generated
by all  the reflections with respect to hyperplanes (in $\A$) containing $F.$ 
So, in general $\W_F$ is not finite. 
\item\label{p3} $U^0:= int(U)$ \ is open in $\R^n$ and a (relatively open) facet
$F\subset\overline{C}_0$ is contained in $U^0$ \ iff \ the stabilizer $\W_F$ is
finite.
\end{enumerate}

 By property \ref{p3} the arrangement  is locally finite in the interior part $U^0.$ 

Let \ $H_{\C}\ :=\ H+iH\ \subset\ \C^n$ \ be the {\it complexification} of the
hyperplane $H,$ and set
$$\mbf{Y}:=\ [U^0\ +\ i\R^n]\ \setminus\ \cup_{H\in\A}\ H_{\C}$$
which corresponds to complexifying only the interior part of the Tits cone.  The
group $\W$ acts (as before) diagonally onto $\mbf{Y},$ and one shows easily 
(by property 3 above) that the action is free. 
Therefore, one has an {\it orbit space} 
$$\mbf{Y_W}:=\ \mbf{Y/W}$$
which is still a manifold, and a regular covering 
$$f_{\W}:\ \mbf{Y}\ \to\  \mbf{Y_W}$$
with group $\W.$

\begin{df}
Define 
$$\GW:=\ \pi_1(\mbf{Y}_W)$$
as the \emph{ Artin group} of type $\mbf{W}.$ 
\end{df}
(see \emph{\cite{brieskorn, deligne}} for the case 
when $\W$ is finite).
We recall here some topological constructions 
from \cite{salvettiArtin} (see also 
\cite{salvetti}, \cite{calmorsal2}, \cite{calmorsal1}, \cite{calmorsal3}, \cite{salvettipreprint}, \cite{morsalvil}).

Take $x_0\in C_0$ and let $Q$ be the {\it finite} $CW$-complex constructed as
follows.

\noindent For all maximal subsets $\Ga\subset S$  such that $\WG$ is finite, construct a
$|\Ga|$-cell  $Q_{\Ga}$ in $U^0$ as the convex hull of the $\WG$ orbit of $x_0.$
Each $Q_{\Ga}$ is a finite convex polyhedron which contains the point $x_0.$

Define 
\begin{equation}\label{poly}
Q\ :=\ \cup\ Q_{\Ga}
\end{equation}
(a finite union of convex polyhedra).


The $k$-faces of each $Q_{\Gamma}$ are also polyhedra, each of them 
corresponding to a coset
of a {\it parabolic} subgroup $\W_{\Delta},$ where $\Delta\subset S$ is a $k$-subset of
$\Ga.$ The correspondence
$$\{ \mbox{faces of $Q_{\Ga}$}\} \ \leftrightarrow\ \{ w.\W_{\Delta},\ w\in\Ga,\ \Delta\subset \Ga\}$$
is obtained by taking the  polyhedron given by the convex-hull of the orbit
$\W_{\Delta}.x_0$ and translating it by 
$w.$

One also has:

\begin{prop}\label{prop:minimal} Inside each coset $w.\W_{\Delta}$ there exists an
unique element  of {\rm minimal length}.
\end{prop}

Here the length is the minimal number of letters (coming from $S$) in a reduced
expression.

For every face $e$ of $Q_{\Ga},$ which corresponds to a coset $w.\W_{\Delta},$ let
$\beta(e)\in w.\W_{\Delta}$ be the element of minimal length.  Notice that $\W_{\Ga}$ 
permutes faces of the same dimension. 
Let $\X_{\W_{\Ga}}$ be obtained from  $Q_{\Ga}$ by identifications on its faces
defined as:  each pair of faces $e,\ e'$ belonging to the same orbit is identified by
using the homeomorphism
$\beta(e)\beta(e')^{-1}.$
Define also the finite complex
\begin{equation}\label{salvcomp}
\mbf{X_W}\ :=\ \cup_{\Ga}\ \X_{\W_{\Ga}}.
\end{equation}

Notice that  $\X_{\W}$  is well defined 
because of the following easy fact.

\begin{rmk}\label{rmk:compatibility} 
For any common cell $e\subset Q_{\Ga}\cap Q_{\Ga'}$ the minimal element 
$\beta(e)$ is the same when computed in $\WG$ and in $\W_{\Ga'}$.
\end{rmk}

For reader convenience, we also recall the following fact, whose proof is the same as in 
\cite{salvettiArtin}.
 
 \begin{teo}\label{teo:main1} The $CW$-complex $\X_{\W}$ is deformation retract
of the orbit space $\mbf{Y_W}.$
 \end{teo}

\begin{dm} First, there exists a regular $CW$-complex $\X\subset \mbf{Y}$ which is
deformation retract of $\mbf{Y},$ and $\X$ is constructed as in \cite{salvetti}.
The construction used there 
(where the starting point was an affine arrangement of hyperplanes) 
works in general because one reduces to the finite case around faces with finite
stabilizer. 
 
The construction of $\X$ can be chosen invariantly with respect to the action of
$\W,$ which permutes cells of the same dimension.  
 
The action on $\X$ being free,  we look at the orbit space $\mbf{X/W}.$ By
remark  \ref{rmk:compatibility}, this reduces to finite cases. \qed\end{dm}


As an immediate corollary of this construction one has a presentation of the
group $\GW$ which generalizes the finite type case (see  also  \cite{vanderlek})

\begin{equation}\label{presentationGW}
\GW\ =\ <g_s,\ s\in S\ |\ g_sg_{s'}g_s\dots = g_{s'}g_{s}g_{s'}\dots> 
\end{equation}
(same number   $m(s,s')$ of factors on each side) where we have to consider 
only pairs such that $m(s,s')$ is finite.

\subsection{Topological constructions for Coxeter groups} 
We refer here essentially to \cite{salvdec2}. 

Consider the subspace arrangement
in $\R^{nd}\ \cong (\R^n)^d$   given by
$$\A^{(d)}:=\ \{ H^{(d)}\ \}$$
where $H^{(d)}$ is the codimensional-$d$ subspace given by
``$d$-complexification'' of the hyperplane $H\in \A:$
$$H^{(d)}:=\ \{ (X_1,\dots,X_d)\ :\ X_i\in\R^n,\ X_i\in\ H\ \}.$$
If $U^0$ is the interior of the Tits cone, as above, then we consider
the space
$$\mathcal U_W^{(d)}:= U_0\times\R^n\times\dots\times\R^n$$
($d-1$ factors equal to $\R^n$), and the configuration space 
\begin{equation}\label{confspace}
\Y:= \ \mathcal U_W^{(d)} \setminus \cup_{H\in\A}\ H^{(d)}.
\end{equation}

As before, the group $\W$ acts freely on $\Y$ and we consider the orbit
space
$$\Y_{\W}\ :=\ \Y/\W.$$
One has:
\begin{teo}\label{kpi1} The space
$$\mbf{Y}_{\mbf{W}}^{(\infty)}:=\ \left[ \lim_{\stackrel{\longrightarrow}{d}}\ \Y\right]
/ \W\ =\  \left[ \lim_{\stackrel{\longrightarrow}{d}}\ \Y_{\W}\right] $$
is a space of type $k(\W,1).$
\qed\end{teo}

So, different from the case of Artin groups, we always get a $k(\pi,1)$ space
here.

Recall also: 
\begin{teo}\label{S-complex}  The space $\Y_{\W}$ contracts
over a $CW$-complex $\X^{(d)}_{\W}$
such that
$$
\{ k-\mbox{cells of } \X^{(d)}_{\W}\}\ \longleftrightarrow \{ \mbox{flags }
\mbf{\Gamma}:=(\Ga_1\supset\dots \supset\Ga_d):\ \Ga_1\subset S, \ \W_{\Ga_1}\mbox{ finite, } \ $$ $$ 
\sum_{i= 1}^d |\Ga_i|\ =\ k\}. $$
Passing to the limit, $\mbf{Y}^{(\infty)}_{\W}=k(W,1)$ contracts over a
$CW$-complex $\X^{(\infty)}_{\W}$
such that
$$
\begin{array}{ccr}
\{ k-\mbox{cells of } \X^{(\infty)}_{\W}\} & \longleftrightarrow & \{ \mbox{flags
} \mbf{\Gamma}:=(\Ga_1\supset\Ga_2\supset\dots):\ \Ga_1\subset S, \\
&& \W_{\Ga_1}\mbox{ finite, } \sum_{i\geq 1} |\Ga_i|\ =\ k\}.
\end{array}
$$

\qed\end{teo}

Notice that $\X^{(\infty)}_{\W}$ does not have finite dimension but the number
of  $k$-cells is finite, given by $\binom {n+k-1}{k}.$

\subsection{Algebraic complexes for Artin groups}\label{subsec:Artin}
We refer here mainly to \cite{salvettiArtin},\cite{salvDecArtin}.

We consider the algebraic complex related to the cell structure of $\X_{\W}.$ It
is given in the following way.

Let $\Z[\GW]$ be the group algebra of $\GW.$  Let $(C_*,\partial_*)$ be the
algebraic complex of free $\Z[\GW]$-modules such that in degree $k$ it is free
with basis $e_J$ corresponding to subsets $J\subset S$ such that $\W_J$ is
finite:

\begin{equation}\label{algcompl}
C_k\ :=\ \bigoplus_{\scriptsize \begin{array}{l} J\subset S \\ |J|=k \\ \W_J \mbox{ finite} 
\end{array}}\ \Z[\GW]\ e_J
\end{equation}

Let

\begin{equation}\label{bordoArtin}
\partial(e_J)\  :=\ \sum_{\scriptsize \begin{array}{c} I\subset J \\ |I|=k-1\end{array}}\
[I:J] \ T^J_I\ . e_I
\end{equation}
where $[I:J]$ is the incidence number ($=0,\ 1$ or $-1$) of the cells in
$\X_{\W}$ and
$$T^J_I\ :=\ \sum_{\beta\in W^J_I}\ (-1)^{\ell(\beta)}\ g_{\beta}$$
where
\begin{enumerate}
\item $\W^J_I\ :=\ \{\beta\in \W_J\ :\ \ell(\beta s)\ >\ \ell (\beta),\
\forall s\ \in\ \W_I \}$  \ 
 is the set of elements of minimal length for  the cosets $\frac{\W_J}{\W_I}$ \
(prop. \ref{prop:minimal});
\item\label{length} if $\beta \in \W^J_I$ and $\beta=s_{i_1}\dots s_{i_k}$ is a
reduced expression then $\ell(\beta)=k;$ 
\item  if $\beta$ is as in \ref{length}  then  $g_{\beta}:= g_{s_{i_1}}\dots
g_{s_{i_k}}.$  One shows that this map 
\begin{equation}\label{section}
\psi: \W \to \GW
\end{equation}
is a well-defined section (not a homomorphism) of the standard surjection 
$\GW \to \W.$ 
\end{enumerate}

\begin{rmk}\label{rmk: risoluzione} When $\XW$ is a space of type $k(\GW,1)$ then $(C_*,\partial_*)$ 
gives a free $\Z[\GW]$ resolution of $\Z$ as a trivial module. In any case, $(C_*,\partial_*)$
corresponds to the cellular complex structure of the universal covering of $\XW,$ 
so they can be used to compute local systems over $\XW.$
\end{rmk}

Let $R:=A[q,q^{-1}]$ be the ring of Laurent polynomials over a ring $A.$ One can
represent $\GW$ by
\begin{equation}\label{rappr}
g_s\ \mapsto \  [\mbox{multiplication by } -q]  \quad \forall s\in S
\end{equation}
($\in Aut(R)$).  

The tensor product $C_*\otimes R$ has boundary
\begin{equation}\label{bordoq}
\partial(e_J)\ =\ \sum_{\scriptsize \begin{array}{c}I\subset J\\ |I|=|J|-1\end{array}}\ [I:J]\ \frac{\W_J(q)}{\W_I(q)}\ e_I
\end{equation}
where
$$\W_J(q):=\ \sum_{w\in W_J}\ q^{\ell(w)}$$
is the Poincar\'e series of the group $\W_J$ (here, a polynomial since the
stabilizers are finite). The denominator $\mathbf W_I(q)$ divides the numerator $\mathbf W_J(q)$
so the quotient is still a polynomial.

\subsection{Algebraic complexes for Coxeter groups}\label{Coxeter}
Consider the algebraic complex
$(C_*, \de)$ of free  
$\Z[\W]$-modules, where
$$
C_k:= \bigoplus_{\scriptsize \begin{array}{c} \mbf{\Gamma}:
\sum_{i\geq 1} \card{\Gamma_i}=k \\ \card{W_{\Gamma_1}}< \infty
\end{array}} \Z[\W] e(\mbf{\Gamma})
$$
The generators of $C_*$ are in one to
one correspondence with the cells of $\X^{(\infty)}_{\W}$, so with the flags
$\mbf{\Ga}= (\Ga_1\supset\Ga_2\supset\dots),\ \Ga_1\subset S, \ \W_{\Ga_1}\mbox{ finite. }$

The
expression of the boundary is the following:
\begin{equation}\label{eq:bordoCoxeter}
\de e(\mbf{\Gamma})= \sum_{ \scriptsize
                            \begin{array}{c}
                                         i \geq 1 \\
                                    \card{\Gamma_i}>\card{\Gamma_{i+1}}
                            \end{array}
                            }
                    \sum_{\tau \in \Gamma_i}
                    \sum_{ \scriptsize
                            \begin{array}{c}
                                    \beta \in \W_{\Gamma_i \setminus \{ \tau
                                    \}}^{\Gamma_i}\\
                                    \beta^{-1} \Gamma_{i+1} \beta
                                    \sst \Gamma_i \setminus \{ \tau
                                    \}
                            \end{array}
                            }
                            (-1)^{\alpha(\mbf{\Gamma}, i, \tau,
                            \beta)}\beta e(\mbf{\Gamma'})
\end{equation}
where
$$
\mbf{\Gamma'}= (\Gamma_1\supset \ldots\supset \Gamma_{i-1}\supset \Gamma_i \setminus
\{\tau\} \supset \beta^{-1} \Gamma_{i+1} \beta\supset \beta^{-1}\Gamma_{i+2} \beta\supset\dots )
$$
and $(-1)^{\alpha(\mbf{\Gamma}, i, \tau, \beta)}$ is an incidence
index. To get a precise expression for ${\alpha(\mbf{\Gamma},i,
\tau, \beta)}$, fix a linear order on $S$ and let
\begin{align*}
\mu(\Gamma_i, \tau):=&\card{j \in \Gamma \textrm{ s.t. } j \leq \tau}\\
\sigma(\beta, \Gamma_j):=& \card{(a,b) \in \Gamma_j\times \Gamma_j
\textrm{ s.t. } a<b \textrm{ and } \beta(a)>\beta(b)}
\end{align*}
in other words, $\mu(\Gamma_i, \tau)$ is the number of reflections
in $\Gamma_i$ less or equal to $\tau$ and $\sigma(\beta,
\Gamma_j)$ is the number of inversions operated by $\beta$ on
$\Gamma_j$. Then we define:

\begin{equation}\label{eq:signbordoCox}
    {\alpha(\mbf{\Gamma}, i, \tau, \beta)}= i \ell(\beta) +
    \sum_{j=1}^{i-1} \card{\Gamma_j} + \mu(\Gamma_i, \tau) +
    \sum_{j=i+1}^{d} \sigma(\beta, \Gamma_j)
\end{equation}
where $\ell$ is the length function in the Coxeter group.

\begin{teo}\label{algcompCox}
For any finitely generated $\W,$ the algebraic complex $(C_*,\de_*)$  gives a
free resolution of the trivial 
$\Z[\W]$-module $\Z.$
\end{teo}
The proof follows straightforward from the remark that the limit space
$\mbf{Y}^{(\infty)}_{\mbf W},$  so 
$\mbf{X}^{(\infty)}_{\mbf W},$
is a space of type $k(\W,1).$

\section{The genus problem for Artin groups}\label{genus}
Our main application here is the extension of some of the results found in
\cite{salvdec2}, \cite{salvdecproc3} about the 
genus of the covering associated to an Artin group.

\subsection{Schwarz genera and homological genera}
We start recalling the definition of Schwarz genus and discussing
 briefly some of its properties (we refer to \cite{schwarz}, \cite{vassiliev} 
 for details).
\begin{df}
For a locally trivial fibration $f: \mbf{Y} \to \X$, the Schwarz genus
$g(f)$ is the minimal cardinality of an open cover $\mathcal{U}$
of $\X$ such that $f$ admits a section over each set $U \in
\mathcal{U}$.
\end{df}

 \noindent {\bf{Remark.}}\ The Schwarz genus is
the extension to fibrations of the Lusternik-Schnirelmann category
of a topological space; indeed the category of a path connected
topological space coincides with the Schwarz genus of its Serre
fibration.\vs

 When $\X$ has the homotopy type of  a finite dimensional CW
 complex, we have an upper bound for the genus of any
 fibration:
\begin{teo}\label{teo:up:genus} If $X$ has the homotopy type of
a CW
complex of dimension $N$, then $g(f) \leq N+1$.
\end{teo}
\qed

 Let now $f: \mbf{Y} \to \X$ be a regular
$G$-covering. Then  $f=a^*(p),$ where  $a: \X \to BG$
is a classifying map into the classifying space  $BG$ for $G$ 
and $p:\ EG \to BG$ is the universal $G$-bundle . 

Let $M$ be an arbitrary $G$-module and $a^* M$ be the local system
on $X$ induced by the map $a$.
\begin{df}The homological $M$-genus of $f:\mbf{Y}\to \X$ is the smallest
integer $h_M(f)$ such that the induced map in cohomology:
$$
a^*:\, H^j(BG; M)\to H^j (\X; a^*M)
$$
is zero in degree $j$ for $j\geq h_M(f)$.\\
The homological genus is defined as the maximum $h(f)=\max_M
h_M(f)$ of the homological $M$-genera.
\end{df}

Homological genus provides a lower bound for Schwarz genus:
\begin{teo}\label{teo:low:genus}
For any regular covering $f: \mbf{Y} \to \X$, we have $g(f) \geq h(f)$.
\qed
\end{teo}

\subsection{The genus problem for Artin groups}
Let $\W$ be a Coxeter group and consider the regular covering $f_{\W}:
\mbf{Y} \to \mbf{Y}_{\W}$ between the configuration spaces introduced in part 
\ref{subsec:topconArt}. We are interested in the genus $g(f_{\W})$ of
$f_{\W}$.  

We start by some general remarks on the algebraic complexes for $\GW$.

First we denote by $K:=K(\W)$ the simplicial scheme, 
defined over $S,$ of the subsets $J\subset S$ which generate 
a finite parabolic subgroup $\W_{\Ga}$ (we include the empty set which by definition generates the trivial subgroup). The algebraic complex which computes the simplicial homology of $K$
with coefficients in $\Z$ will be denoted here by

$$D_k^0(\W):= \bigoplus_{\scriptsize \begin{array}{c} J \subseteq  S\\
 \card{J}=k
\\ \card{W_{J}}< \infty
\end{array}} \Z \cdot e_J^0
$$
with boundary:
\begin{align}\label{eq:simplicial:boundary}
\de^0 (e_J^0)
=&\sum_{\scriptsize \begin{array}{c} I \subset J \\
\card{J}=\card{I}+1
\end{array}}            [I:J]\,
 e_I^0
\end{align}

Here we indicate by $[I:J]$ the incidence number of oriented simplices, namely $\pm 1$ or $0,$ which is the same as that of the corresponding cells appearing in  (\ref{bordoq}). 

Notice that we found convenient here to graduate $D^0_*$ 
according to the cardinality of the subsets, so there is a degree$-1$ shift isomorphism with the standard complex for the simplicial homology of $K:$ \ $H_m(D^0_*) \ \cong\ H_{m-1}(K),\quad m\geq 1.$

We also denote by $D^0_*(\W,B):=\ D^0_*(\W)\otimes B$ the algebraic complex computing the homology with trivial coefficients a $\Z-$module $B.$ When no module is indicated we mean $B=R.$

Let us consider the representation  $\rho: \GW\to Aut(R)$ of (\ref{subsec:Artin}), 
obtained by sending the standard generators of $\br$ into
$(-q)$-multiplication. Let $R_q$ be the the ring $R$
with the prescribed structure of $\br$-module. 
It is convenient to indicate here by $D_*(\W)$ the algebraic complex 
of part \ref{subsec:Artin}, so
\begin{equation}\label{eq:chain:compl:rank1}
D_k(\W):= \bigoplus_{\scriptsize 
\begin{array}{c} J \subseteq  S\\
 \card{J}=k
\\ \card{\W_{J}}< \infty
\end{array}} R \cdot e_J
\end{equation}
and the boundary is given by (\ref{bordoq}).

\begin{rmk}  We can formally
rewrite the boundary map in (\ref{bordoq}) as:
\begin{align}\label{eq:formal:boundary}
\de (\frac{1}{\W_J(q)}e_J)
=&\sum_{\scriptsize \begin{array}{c} I \subset J \\
\card{J}=\card{I}+1
\end{array}}            [I:J]
\frac{1}{\W_I(q)} \cdot e_I
\end{align}
\end{rmk}
That means that the fractions \ ${e_J}/{\W_{J}(q)}$\ 
behave like the cells of the simplicial scheme  $K.$\\
Consider the diagonal map:
\begin{align}
\Delta: D_* (\W)&\to D_*^0(\W), \quad  e_J  \mapsto \W_J(q) e_J^0.
\end{align}
It is clear by the previous discussion that $\Delta$ is an injective
chain-complex homomorphism, so  
there is an exact sequence of complexes:
\begin{equation}\label{eq:short:rank1}
\xymatrix @R=.5pc @C=3pc {
    0 \ar[r] & D_* (\W)\ar[r]^-{\Delta} & D_*^0 (\W)\ar[r]^-{\pi} & L_*(\W)
    \ar[r] & 0
}
\end{equation}
where
$$L_k(\W):= \bigoplus_{\scriptsize \begin{array}{c} J \subseteq  S\\
 \card{J}=k
\\ \card{\W_{J}}< \infty
\end{array}} \frac{R}{(\W_J(q))} \cdot \bar{e}_J
$$
is the quotient complex.

Passing to the associated long exact sequence we get:
\begin{equation}\label{eq:long:rank1}
\xymatrix @R=.5pc @C=1.5pc {
    \ar[r]^-{\pi_*}&H_{k+1}(L_*) \ar[r]^-{\delta} & H_k(D_*) \ar[r]^-{\Delta_*}
& H_k(D_*^0) \ar[r]^-{\pi_*}&
    H_k(L_*)
    \ar[r]^-{\delta} & H_{k-1}(D_*) \ar[r]^-{\Delta_*}&
}
\end{equation}

\begin{rmk}\label{weightedsheaves}
It is possible to consider both  the data  $\frac{R}{(\W_J(q))},$ $J\subset S, \ |\W_J|<\infty,$
and the complex $L_*$ as  \emph{ functorial constructions} associated to any (finitely generated) 
Artin group. In \cite{morsalvil} we introduce a class of ``sheaves over posets" called \emph{ weighted
sheaves over posets}, and associated \emph{ weighted complexes}, a particular case being that associated to an Artin group,  and we used this construction  for computations of the cohomology.
\end{rmk}

We need some definitions.

\begin{df} We define the \emph{virtual dimension} of an Artin group $\GW$ as 
$$vd(\GW)\ :=\ dim (K)+1$$
where $K$ is the associated simplicial scheme and the dimension is that of a simplicial complex,
so  by definition of $K:$
$$vd(\GW)\ =\ max\{ |J|\ :\ J\subset S, \ |\W_J|<\infty \}$$
\end{df} 

\bigskip

\noindent (Equivalently, $vd(\GW)=\ max\{n\ :\ D^0_n(\W;\Z)\neq 0\}$)
\bigskip

\begin{df} We define the \emph{homological virtual dimension} of $\GW$ as:
$$hvd(\GW)\ :=\  max\{n:\ H_n(D^0_*(\W;\Z))\neq 0\ \}$$
\end{df}
\bigskip

\noindent (Equivalently: $hvd(\GW)\ :=\ max\{n\ :\ H_n(K;\Z)\neq0\}\ +1$,\ where here 
we use the standard graduation for the homology).

\bigskip
By construction one has
$$dim(\XW)\ =\ vd(\GW)$$
and from theorem \ref{teo:up:genus} it follows

\begin{equation}\label{equation:upbound}
g(f_{\W})\ \leq \ vd(\GW)+1
\end{equation}
 
\bigskip

\begin{df} We say that  a Coxeter system $(\W,S)$ is  \emph{affine-like} type if
$$vd(\GW)\ =\ hvd(\GW)$$
\end{df}
 
\bigskip

\begin{rmk}\label{rmk:affine} \ Recall that for an affine Coxeter
system $(\W,S)$ of rank $n+1$, a parabolic subgroup $\W_J$ is finite
if and only if $J$ is a proper subset of $S$. In particular the
poset of finite parabolic subgroups is isomorphic to the poset of
proper subsets of $I_{n+1}=\{1, \ldots, n+1\}$ (that is the
boolean lattice minus its maximum). 
Then the homology of $D^0_*(\W)$ is the
reduced homology of a $(n-1)$-sphere modulo a degree shift:
$$
H_k(D^0_*(\W))\cong \t{H}_{k-1}(S^{n-1}; R)\cong \left\{
\begin{array}{cl} 0 & \textrm{if $k\neq n$}\\ R &\textrm{if $k=n$} \end{array} 
\right.
$$
Therefore $(\W,S)$ is affine-like.
\end{rmk}

\bigskip
We can now state the main result of the paper. 	
\bigskip

\begin{teo}\label{teo:main} Let $(\W,S)$ be an affine-like Coxeter system. Then for the Schwarz genus of the fibration 
$$f_{\W}:\ \mbf{Y} \to \mbf{Y}_{\W}$$
it holds
$$g(f_{\W})\ =\ vd(\GW)+1$$
\end{teo} 

\begin{dm}  Inequality (\ref{equation:upbound}) gives the upper bound which we need. To obtain
 the lower bound we will use homological methods.
 
First we proof:

\begin{teo}\label{thm:rank:affine}
Let  $M=\Z[-1]$ be the $\GW$-module $\Z$ with the action given by the sign 
representation.  Let   
$$F^0_k:=\ FH_k(D^0_*(\W;\Z))\ :=\ H_k(D^0_*(\W;\Z))/Tors(H_k(D^0_*(\W;\Z))$$  
be the free component of the integral $k-$th homology of $D^0_*.$ 
Assume that 
$$F^0_k\ \neq\ 0$$ 
for some $k.$ 
 
Then also the free part
$$FH_k(\X_W;M)\ \neq\ 0.$$ 
\end{teo}
\begin{dm}[Proof of teorem \ref{thm:rank:affine}]  The homology of $\XW$ in 
the sign representation is computed by specializing the complex 
(\ref{eq:chain:compl:rank1}) to $q=1.$  We obtain sequences analog to 
(\ref{eq:short:rank1}) and (\ref{eq:long:rank1}) respectively, where in this case 
$D^0_*(\W;\Z)$ computes the homology of $K$ with trivial integer coefficients.

By definition $F^0_k\neq 0.$ Let $z^0\in D^0_k$ be any no-torsion $k$-cycle, 
$$z^0\ =\ \sum_{|J|=k}\ \epsilon_J\ e^0_J,\ \quad \epsilon_J\in\Z.$$

Notice that the map $\Delta$ is injective, being diagonal with all non-vanishing entries
given by $\W_J(1)\neq 0,$ $|J|=k.$

Let $\mu:=\ \mbox{lcm }(\W_J(1))_{|J|=k}.$  Then

$$\mu.z^0\ =\ \Delta(z)$$
where
$$z:=\ \sum_{|J|=k}\  \epsilon_J \frac{\mu}{\W_J(1)}\ e_J\ .$$

It follows that $z$ is a no-torsion cycle in $H_n(\XW;M)$  
 \qed \end{dm} 

\bigskip

From theorem \ref{thm:rank:affine} we deduce:

\begin{cor}\label{cor:rank:affine}  Let $(\W,S)$ be an affine-like Coxeter system and let 
$$n:= \ vd(\GW)\ =\ hcd(\GW).$$
Let $M=\Z[-1]$ be the sign representation. 

Then the $\Z$-rank of $H^n(\XW;M)$ is bigger than $0.$
\end{cor}
\begin{dm}[Proof of corollary \ref{cor:rank:affine}.]  
 In top dimension $n$ we have (we drop the coefficients):

\begin{equation}\label{eq:longsign:rank1}
\xymatrix @R=.5pc @C=1.5pc {
    0 \ar[r]  & H_n(D_*) \ar[r]^-{\Delta_*}
& H_n(D_*^0) \ar[r]^-{\pi_*}&
    H_n(L_*)
    \ar[r]^-{\delta} & H_{k-1}(D_*) \ar[r]^-{\Delta_*}&
}
\end{equation}

By definition $H_n(D^0_*)\neq 0,$ and since we are in top dimension such group is
a free $\Z$-module. Then we apply theorem \ref{thm:rank:affine} and
the thesis follows by passing to cohomology.
 \qed \end{dm}

Recall that we have an inclusion $i:\X_{\W} \hookrightarrow
\X^{\pp{\infty}}_{\W}$ and that $\X_{\W}$ may be identified with the
subcomplex of $\X^{\pp{\infty}}_{\W}$ consisting of cells of type
$\mbf{\Gamma}=(\Gamma_1\tss \emptyset \tss \emptyset \tss
\ldots)$.

Let $M$ be a $\W$-module and $M'$ the local coefficient system on
$\X_{\W}$ induced by $M$ via $i$.  Equivalently, we can consider the natural surjection
$$\GW \to \W\to 1.$$
The associated map of cochains
$$
i^*: C^*({\W}; M) \to C^*(\X_{\W}; M')
$$
is described as the restriction of $c \in
C^*({\W};M)$ to the chains for $\X_{\W}$. Let $n$ be the maximal
cardinality of a subset $J\subset S$ s.t. $\card{{\W}_J}< \infty$.
Then in degree $n$ we have:
\begin{equation}\label{eq:comm:diagram:top}
\xymatrix@C=1.7em{\ar[r]& C^{n-1}({\W}; M) \ar[r]\ar[d] & C^n({\W}; M)
\ar[r]\ar[d] & C^{n+1} ({\W};
M)\ar[r]\ar[d] &\\
\ar[r]& C^{n-1}(\X({\W}); M') \ar[r] & C^n(\X_{\W}; M') \ar[r] & 0\ar[r]
&}
\end{equation}

The above sign representation $\Z[-1]$ is clearly induced by the analog sign representation of $\W.$
We have:

\begin{teo}\label{teo:epimo:top:sign}
Let $M=\Z[-1]$ be the sign representation as above. Then the
map $i^*: H^k({\W}; M) \to H^k(\X_{\W}; M')$ is an epimorphism, for all $k\geq0.$
\end{teo}

\begin{dm}[Proof of theorem \ref{teo:epimo:top:sign}.]  
The proof is entirely analog to \cite{salvdec2}, where only the case $k=n$ is considered.
We include it here for reader convenience. 
 
We consider the diagram \ref{eq:comm:diagram:top} at vertical level $k.$

Let $S^{(\infty)}_k$ be the set of (infinite) flags of total cardinality $k$ (see part \ref{Coxeter})
 so that we identify 
$C^k(\W;M)$ with the set of functions  $f: S^{(\infty)}_k\to M,$ denoted
 $M^{S^{(\infty)}_k}.$  Analogously we identify  $C^k(\XW;M')$ with the set 
 of such functions defined only over the flags of length $1,$  i.e. such that $\Ga_i=\emptyset$ for $i>1.$ 
 
 It is sufficient to consider functions \ $f\in  M^{S^{(\infty)}_n}$ \ such that

$$f(\mbf{\G}=(\G_1\supset\G_2\supset\dots))\ =\ 0 \quad \text{if} \ |\G_i|>0 \
\text{for some}\ i>1$$

\ni (equivalently, if $|\G_2|>0$). Let us compute the coboundary of $f$.

\begin{enumerate}

\item[(a)]   $\delta^k(f)(\G_1\supset\G_2\supset\dots)\ =\ 0 $  \quad  if \
$|\G_2|>1$ \ or \ $|\G_3|>0$.

In fact, by formula (\ref{eq:bordoCoxeter})  we get a linear combination of $f$ computed on chains
with $|\G_2|>0$.

\item[(b)] If \ $\G_1 = \{ s_{j_1}<\dots <s_{j_k}\},\ \G_2 = \{ s_{j_m} \} $ \
then

 $\delta^k(f) ( \G_1\supset\G_2\supset\emptyset\supset\dots )\ =$

$$=\ \sum_{\beta_2\in\W_{s_{j_m}}} (-1)^{\alpha(\GF,2,1,\beta_2)}\
\rho(\beta_2).f(\G_1)\ =$$

$$=\ (-1)^{k+1}\ f(\G_1)\ +\ (-1)^{k+1}\rho(s_{j_k}).f(\G_1)\ =$$

\begin{equation}\label{eq:bordosign} 
=\ (-1)^{k+1}\ (1+\rho(s_{j_m})).f(\G_1)
\end{equation}

\end{enumerate}
It immediately follows from (\ref{eq:bordosign}) that each cocycle $f\in C^k(\XW;M')$ can be extended to a cocycle in $C^k(\W;M),$ which concludes the proof.
\qed \end{dm}

\bigskip

\begin{dm}[End of proof of theorem \ref{teo:main}.]
From Theorem
\ref{thm:rank:affine}, we know that  the top-cohomology of
$\X_{{\W}}$ with coefficients in the sign representation does not
vanish. Using theorem \ref{teo:epimo:top:sign}, the
homological genus $h(f_{{\W}})$ is greater than $vd(\GW)+1$. Since
$\X_{{\W}}$ has dimension $vd(\GW)$, the result follows from
theorems \ref{teo:low:genus}
and \ref{teo:up:genus}.
\qed
\end{dm}
\end{dm}

\bigskip
We have immediately (see rmk \ref{rmk:affine}):

\begin{cor}\label{teo:genus}
Let ${\W}_a$ be an affine Weyl group of rank $n+1$. Then the Schwarz
genus of the fibration $\mbf{Y}_{{\W}_a} \to \X_{{\W}_a}$ is precisely $n+1$.
\qed
\end{cor}
\bigskip

In case  $hvd(\GW)$ is strictly lower than $vd(\GW)$ we can get an estimate by 
slightly modifying the definitions.

\bigskip

\begin{df} We define the \emph{rational homological virtual dimension} of $\GW,$ written \ 
$rhvd(\GW),$ as
$$rhvd(\GW)\ :=\ max \{n:\ H_n(D^0_*(\W;\mathbb Q))\neq 0\ \}.$$
\end{df}

\bigskip

Then a straightforward modification of the proof of theorem \ref{teo:main}  implies:
\bigskip

\begin{teo} For any finitely generated Coxeter system $(\W,S)$ one has 
$$g(f_{\W})\ \geq\  rhvd(\GW)+1.$$
\qed\end{teo}


\providecommand{\bysame}{\leavevmode\hbox
to3em{\hrulefill}\thinspace}
\providecommand{\MR}{\relax\ifhmode\unskip\space\fi MR }
\providecommand{\MRhref}[2]{%
  \href{http://www.ams.org/mathscinet-getitem?mr=#1}{#2}
} \providecommand{\href}[2]{#2}


\begin{thebibliography}{DPSS99}


\bibitem[Aro05]{arone}
G.~Arone, \emph{A note on the homology of {$\Sigma_n$}, the {S}chwartz
  genus, and solving polynomial equations}, Proceedings of the Arolla Conference on Algebraic Topology, (2005).


\bibitem[Bou68]{bourbaki}
N.~Bourbaki, \emph{Groupes et algebr\`es de {L}ie}, vol. Chapters {IV-VI},
  Hermann, (1968).


\bibitem[Bri71]{brieskorn}
E.~Brieskorn, \emph{Die {F}undamentalgruppe des {R}aumes der regul{\"a}ren
  {O}rbits einer endlichen komplexen {S}piegelunsgruppe}, Invent. Math.
  \textbf{12} (1971), 57--61.




\bibitem[CMS08a]{calmorsal2}
F.~Callegaro, D.~Moroni, and M.~Salvetti, \emph{Cohomology of affine {A}rtin
  groups and applications}, TRANSACTIONS OF THE AMERICAN MATHEMATICAL SOCIETY
  \textbf{360} (2008), 4169--4188.

\bibitem[CMS08b]{calmorsal1}
\bysame, \emph{Cohomology of artin groups of type $\tilde{A}_n, b_n$ and
  applications}, GEOMETRY \& TOPOLOGY MONOGRAPHS \textbf{13} (2008), 85--104.

\bibitem[CMS10]{calmorsal3}
\bysame, \emph{The $k(\pi,1)$ problem for the affine artin group of type tilde
  $b_n$ and its cohomology,}, JOURNAL OF THE EUROPEAN MATHEMATICAL SOCIETY
  \textbf{12} (2010), 1--22.

\bibitem[Del70]{deligne}
P.~Deligne, \emph{Equation diff\'erentielles \`a points singuliers
  r\'eguliers}, Lecture Notes in Mathematics, vol. 163, Springer-Verlag, (1970).

\bibitem[DPS04]{salvdecproc3}
C.~{De Concini}, C.~Procesi, and M.~Salvetti, \emph{On the equation of degree $6$}, Comm. Math. Helv. \textbf{79}
  (2004), 605--617.



\bibitem[DS96]{salvDecArtin}
C.~{De Concini} and M.~Salvetti, \emph{Cohomology of {A}rtin groups}, Math.
  Res. Lett. \textbf{3} (1996), 293--297.

\bibitem[DS00]{salvdec2}
\bysame, \emph{Cohomology of {A}rtin groups and {C}oxeter groups}, Math. Res.
  Lett. \textbf{7} (2000), 213--232.





\bibitem[Hum90]{humphreys}
J.E. Humphreys, \emph{Reflection groups and {C}oxeter groups}, Cambridge
  University Press, (1990).



\bibitem[Mor06]{moronitesi}
D.~Moroni, \emph{Finite and infinite type artin groups: Topological aspects and
  cohomological computations}, PhD thesis (2006).


\bibitem[MSV12]{morsalvil}
 D.~Moroni, M.~Salvetti and A. Villa \emph{Some topological  problems on the configuration spaces of Artin and Coxeter groups }, in ``Proceedings
of {\it Configuration spaces: geometry, combinatorics and topology}'', CRM E.De Giorgi, (2010), Edizioni della Scuola Normale Superiore.


\bibitem[Sal87]{salvetti}
M.~Salvetti, \emph{Topology of the complement of real hyperplanes in
  $\mathbb{C}^n$}, Invent. Math. \textbf{88} (1987), no.~3, 603--618.

\bibitem[Sal94]{salvettiArtin}
\bysame, \emph{The homotopy type of {A}rtin groups}, Math. Res. Lett.
  \textbf{1} (1994), 567--577.

\bibitem[Sal05]{salvettipreprint}
\bysame, \emph{On the cohomology and topology of Artin and Coxeter groups},
  Pubblicazioni Dipartimento di Matematica L.Tonelli, Pisa (2005).

\bibitem[Sch61]{schwarz}
A.~S.~Schwarz, \emph{Genus of a fibre bundle}, Trudy Moscow. Math. Obshch
  \textbf{10} (1961), 217--272.

\bibitem[Vas92]{vassiliev}
V.~A. Vassiliev, \emph{Complements of discriminants of smooth maps: Topology and
  applications}, Translations of Mathematical Monographs, vol.~98, AMS, (1992).

\bibitem[vdL83]{vanderlek}
H.~van~der Lek, \emph{The homotopy type of complex hyperplane complements},
  Ph.D. thesis, University of Nijmegan, (1983).

\bibitem[Vin71]{vinberg}
E.B. Vinberg, \emph{Discrete linear groups generated by reflections}, Math.
  USSR Izvestija \textbf{5} (1971), 1083--1119.




\end{thebibliography}
\end{document}